\documentclass[12pt]{amsart}
\usepackage{amsmath,amssymb,amsthm,verbatim,hyperref}

\DeclareMathOperator{\dist}{dist}

\begin{document}
\newtheorem{thm}{Theorem}[section]
\newtheorem{lem}[thm]{Lemma}
\newtheorem{dfn}[thm]{Definition}
\newtheorem{cor}[thm]{Corollary}
\newtheorem{conj}[thm]{Conjecture}
\newtheorem{clm}[thm]{Claim}
\theoremstyle{remark}
\newtheorem{exm}[thm]{Example}
\newtheorem{rem}[thm]{Remark}
\def\N{{\mathbb N}}
\def\G{{\mathbb G}}
\def\Q{{\mathbb Q}}
\def\R{{\mathbb R}}
\def\C{{\mathbb C}}
\def\P{{\mathbb P}}
\def\Z{{\mathbb Z}}
\def\v{{\mathbf v}}
\def\x{{\mathbf x}}
\def\O{{\mathcal O}}
\def\M{{\mathcal M}}
\def\kbar{{\bar{k}}}
\def\tr{\mbox{Tr}}
\def\id{\mbox{id}}
\def\qed{{\tiny $\clubsuit$ \normalsize}}

\renewcommand{\theenumi}{\alph{enumi}}

\title{Rational approximation on cubic hypersurfaces}

\author{David McKinnon}

\begin{abstract}
We compute the constant of approximation for an arbitrary rational point on an arbitrary smooth cubic hypersurface $X$ over a number field $k$, provided that there is a $k$-rational line somewhere on $X$.  In the process, we verify the Coba conjecture for $X$.
\end{abstract}

\maketitle

\vspace{-5pt}
\section{Introduction}

Let $X$ be a smooth projective variety defined over a number field $k$.  In \cite{McK}, the author makes a conjecture that the best rational approximations to a rational point $P$ on a smooth projective variety can be found along some rational curve through $P$.  Such a curve is called a curve of best approximation to $P$.

This is a little vague, obviously, but -- as happens so frequently in mathematics -- we must wade through some tedious definitions to make it precise.  To wit: we must say what we mean by a good approximation.

Generally speaking, an approximation  to $P$ is good (in this sense) if:
\begin{enumerate}
	\item the approximating point is close to $P$, and
	\item the approximating point has small height.
\end{enumerate}
Let's quantify that.

To define distance, we choose a place $v$ of $k$, and use the distance on $X$ induced by $v$.  That is, we embed $X$ in projective space and restrict the standard $v$-adic distance there to $X$.  (It turns out it doesn't matter which embedding we pick.  More on that later.)

To define the height, we choose a divisor class $L$ on $X$ -- probably an ample one -- and choose some (multiplicative, relative to $k$) Weil height $H_L$ associated to $L$.

The tricky part is measuring the balance of these two factors.  To do this, we work in aggregate over a sequence of rational points.  

Let $\{x_i\}$ be a sequence of rational points that converges $v$-adically to $P$.  The $L$-approximation constant of $\{x_i\}$ is $\alpha(\{x_i\},P,L)$, a real number whose definition is tedious enough to put it off to the next section.  Nevertheless, $\alpha(\{x_i\},P,L)$ has the property that the larger it is, the worse the sequence $\{x_i\}$ $v$-adically $L$-approximates $P$.  It is, if you like, the cost of approximating $P$ with the $\{x_i\}$.

We define the $L$-approximation constant $\alpha(P,L)$ to be the infimum of $\alpha(\{x_i\},P,L)$ over all nonconstant sequences $\{x_i\}$.  

Thus, the lower $\alpha(\{x_i\},P,L)$ is, the better the sequence approximates $P$.  The lower $\alpha(P,L)$ is, the better $P$ can be approximated.

Time to deliver on my promise.

\begin{conj}[Coba conjecture]\label{conj:ratcurve}
Let $X$ be a smooth variety defined over a number field $k$, and $P\in X$ a $k$-rational point.  Let $A$ be an ample line bundle on $X$.  If $\alpha(P,A)<\infty$, then there is a curve of best $A$-approximation to $P$.  In other words, there is a curve $C$ such that
\[\alpha(P,A|_C)=\alpha(P,A)\]
where the first $\alpha$ is computed on $C$, and the second is computed on $X$.
\end{conj}

To paraphrase the conjecture still further: The best $A$-approximations to $P$ approximate $P$ just as well as the best $A$-approximations to $P$ that happen to lie on $C$.

This conjecture has been verified in a large number of cases, and Vojta's conjecture is known to imply it in a much larger set of cases.  See \cite{McK} for a discussion of this.  

In forthcoming work of Lehmann, McKinnon, and Satriano, a framework is proposed to deduce the~\nameref{conj:ratcurve} from Vojta's Conjecture for all smooth projective varieties.  This framework uses the machinery of the Minimal Model Program to reduce from a general smooth variety to a Mori fibre space (though not necessarily a smooth one).  Thus, it's natural to look to Mori fibre spaces to verify the~\nameref{conj:ratcurve}.  In this paper, as warned in the title, we consider cubic hypersurfaces.  In dimension three or greater, these are all Mori fibre spaces when they're smooth.

In particular, we verify the~\nameref{conj:ratcurve} for smooth cubic hypersurfaces that contain a $K$-rational line, in Theorem~\ref{cubic}.  It is known by work of Dietmann and Wooley (Theorem~2 of \cite{DW}) that any smooth cubic hypersurface of dimension at least 35 defined over a number field $k$ contains a $k$-rational line.  Work of Brandes and Dietmann (Theorem~1.1 of \cite{BD}) improves this to dimension at least 29 when $k=\Q$, so it is likely that the hypothesis in Theorem~\ref{cubic} is satisfied in a large proportion of cases.

The structure of the paper is quite simple.  Section 2 contains some tedious prerequisites like the definition of $\alpha$.  Section 3 first considers the cases of curves and surfaces separately in Theorem~\ref{cubicsurface}, and then proves the higher-dimensional case in Theorem~\ref{cubic}.  

The author thanks Julia Brandes and Mike Roth for help and advice during the preparation of this paper.

\section{Preliminaries}

The purpose of this section is to make precise our terminology.  Sigh.  Let's get to it.

We start with a projective variety $X$ defined over a number field $k$.  For any line bundle $L$ on $X$, we associate a multiplicative Weil height $H_L$.  Details about how to do this can be found in a number of standard references, such as section 1.1 of \cite{Se}.  There is a small ambiguity in the choice of height, but for any two heights $H_L$ and $H'_L$ constructed in this way, the function $H_L/H'_L$ is a bounded, positive function compactly supported away from zero.  And this ambiguity will not affect the definition of $\alpha$, as we will shortly see below.

Next, we fix a place $v$ of $k$.  We need to define a $v$-adic distance function on $X$.  As hinted in the introduction, we do this by choosing an embedding of $X$ in $\P^n$, and restrict the usual $v$-adic distance on $\P^n$ to $X$.  The details of this are surprisingly unpleasant and are covered in \cite{MR}, section 2.  The important points are (a) defining the $v$-adic distance can be done, and (b) any two distances constructed in this way are equivalent, in that $\dist/\dist'$ is a bounded, positive function compactly supported away from zero.  As with the heights, this will not affect the definition of $\alpha$.

And now for $\alpha$.  Unfortunately, the path to the definition of $\alpha$ is not a pleasingly short one.  We begin by defining the approximation constant for a single sequence.  

\begin{dfn}\label{dfn:seqappconst}
	Let $\{x_i\}$ be a sequence of $k$-rational points that converges $v$-adically to $P$.  We define
	\[A(\{x_i\}, L) = \left\{{
		\gamma\in\R \mid
		\dist(P,x_i)^{\gamma} H_{L}(x_i)\,\,\mbox{is bounded from above}
	}\right\}.
	\]
	The constant of $L$-approximation of $\{x_i\}$ to $P$ is 
	\[\alpha(\{x_i\},P,L)= \inf A(\{x_i\},L)\]
	In particular, if $A(\{x_i\},L)=\emptyset$ then $\alpha(\{x_i\},P,L)=\infty$.
	We call $\alpha(\{x_i\},P,L)$ the $L$-approximation constant of $\{x_i\}$.
\end{dfn}

The idea behind $\alpha$ is that it's the exponent you need on $\dist(x_i,P)$ to make $\dist(x_i,P)^\alpha$ approximately equal to $H_L(x_i)$.  The preceding painful definition is merely making that precise.

Also, as advertised earlier, if we modify either the height or the distance by a bounded, positive function compactly supported away from zero, then the value of $\alpha$ is unchanged.  

That's $\alpha$ for a sequence.  Now we need $\alpha$ for a point, independent of the sequence.  It is, unsurprisingly, the infimum of all the approximation constants of sequences.

\begin{dfn}
	The $L$-approximation constant of $P$ is defined to be
	\[\alpha(P,L)=\inf_{\{x_i\}}\alpha(\{x_i\},P,L)\]
\end{dfn}

See?  Told you.

Next up, a couple of results from \cite{MR2} that will come in handy later.  First, a result about how to compute $\alpha$ on a curve.

\begin{thm}[Theorem~2.6 from \cite{MR2}]\label{thm:curve}
	Let $C$ be any curve birational over $k$ to $\P^1_k$, and $\varphi\colon\P^1\rightarrow C$ the normalization map.
	Then for any ample line bundle $L$ on $C$, and any $P\in C(\overline{\mathbb{Q}})$ we have the equality:
	\[\alpha(P,L|_C)=\min_{q\in \varphi^{-1}(P)} d/r_{q} m_{q}\]
	where $d=\deg(L)$, $m_{q}$ is the multiplicity of the branch of $C$ through $P$ corresponding to $q$, and
	\[r_{q}=
	\begin{cases}
		0 & \text{if $\kappa(q)\not\subseteq k_v$} \\
		1 &\text{if $\kappa(q)=k$} \\
		2 &\text{otherwise.}
	\end{cases}
	\]
	where $\kappa(q)$ denotes the field of definition of $q$.
\end{thm}

The other result is a Liouville-type result that provides lower bounds for $\alpha$.  This is a special case of Theorem 3.3 from \cite{MR2}.  

\begin{thm}\label{thm:Liouville-bound}
	Let $X$ be an algebraic variety defined over $k$, $P\in X(k)$ any $k$-rational point.
	
	Let $\pi\colon\tilde{X}\to X$ be the blowup of $X$ at $P$ with exceptional divisor $E$.
	Let $L$ be a nef line bundle on $X$, and
	$\gamma>0$ a rational number such that $L_{\gamma}:=\pi^{*}L - \gamma E$ is in the effective cone
	of $\tilde{X}$.  Finally let $B'$ be the stable base locus of
	$L_{\gamma}$ and set $B=\pi(B')$.  
	
	\begin{enumerate}
	
		\item There is a positive real constant $M$ such that for all $y\in X(k)-B(k)$, we have $H_{L}(y)\mbox{dist}_v(x,y)^{\gamma}\geq M$.
	
		\medskip
	
		\item For any sequence $\{x_i\}\to P$ of $k$-rational points approximating $P$, if infinitely many points of $\{x_i\}$ are
		outside $B$ then $\alpha(\{x_i\},P,L) \geq \gamma$.
		
		\medskip
		
		\item If $\alpha(P,L) < \gamma$ then $P\in B$ and $\alpha(P,L) = \alpha(P,L|_B)$.
		
		\medskip
		
		\item If $P\in B$ and $\alpha(P,L|_{B})\geq \gamma$ then $\alpha(P,L)\geq \gamma$.
	\end{enumerate}
	
\end{thm}

In other words, Theorem~\ref{thm:Liouville-bound} says that if $L_\gamma$ is effective, then $\alpha(P,L)\geq\gamma$, except maybe for sequences that lie in the image of the stable base locus of $L_\gamma$.  

The idea is to use Theorem~\ref{thm:Liouville-bound} to show that most points don't approximate $P$ very well, and then Theorem~\ref{thm:curve} to show that the terrific curve $C$ we will find actually {\em does} approximate $P$ well.  

\section{Main Theorem}

Let $k$ be a number field, and let $v$ be a place of $k$.  Let $X\subset\P^n$ be a smooth cubic hypersurface defined over $k$, and let $L$ be the hyperplane class on $X$.  

If the dimension of $X$ is one, then of course there cannot be any line on $X$, $k$-rational or not, so the claim is vacuous.  Nevertheless, it is classical\footnote{See for example the second theorem on page 98 of \cite{Se}.} that $\alpha(P,L)=\infty$ for every point $P$ of $X$, and the~\nameref{conj:ratcurve} is vacuously true for $X$.

If the dimension of $X$ is two, then the question of computing $\alpha(P,L)$ is treated in \cite{MR2}.  In \cite{MR2}, the Picard group of the surface is assumed to be entirely defined over $k$, but the argument for computing $\alpha(P,L)$ for the anticanonical class -- that is, the class that embeds $X$ as a cubic hypersurface -- depends only on the existence of a $k$-rational conic through $P$, which is equivalent to the existence of a $k$-rational line on $X$ ... unless $P$ is the intersection of two conjugate irrational lines, in which case the argument breaks down.  We remedy this failing with the following theorem, which also relaxes the smoothness hypothesis on the hypersurface.

\begin{thm}\label{cubicsurface}
Let $X\subset\P^3$ be a cubic surface (irreducible, but possibly singular) defined over $k$.  Let $P\in X(k)$ be any $k$-rational smooth point, and let $L$ be the hyperplane class.  Assume that there is a line $T$ on $X$ defined over $k$, and let $S_P$ be the intersection of $X$ with the tangent plane at $P$.  Assume further that the singular locus of $X$ is supported on $T$.  Then:
\[\alpha(P,L)=\left\{\begin{array}{ll}
1 & \mbox{if $P$ lies on a $k$-rational line} \\
2 & \mbox{if $P$ is an isolated point in $S_P(k)$} \\ 
& \mbox{or if $S_P$ has $k$-rational tangent lines at $P$} \\
3/2 & \mbox{otherwise}
\end{array}\right.\]
In all cases, there is a curve of best approximation to $P$.
\end{thm}

\noindent
{\it Proof:} Let $\pi\colon Y\to X$ be the blowup of $X$ at $P$, with exceptional divisor $E$.  Since $\pi^*L-2E$ is effective with stable base locus $S_P$, Theorem~3.3 of \cite{MR2} immediately shows that $\alpha_L(P)\geq 2$, except possibly for sequences contained in $S_P$.  

Now, $S_P$ is a plane curve of degree $3$ with a singularity at $P$, so we can compute $\alpha(P,L|_{S_P})$ using Theorem~\ref{thm:curve}.  If $P$ is $v$-adically isolated on $S_P$, then of course $\alpha(P,L|_{S_P})=\infty$, and so $S_P$ is no threat to the computation of $\alpha$.

Otherwise, we can assume that $P$ is the $v$-adic limit of $k$-rational points on $S_P$.  Theorem~\ref{thm:curve} confidently lists the possibilities here:

\begin{itemize}
	\item If $S_P$ is reducible or non-reduced, then it's the union of two $k$-rational lines and $\alpha(P,L|_{S_P})=1$.
	\item If $S_P$ is irreducible and there are two tangent lines to $S_P$ at $P$, both defined over $k$, then $\alpha(P,L|_{S_P})=3$.
	\item If $S_P$ is irreducible and there are two tangent lines to $S_P$ at $P$, neither defined over $k$, then $\alpha(P,L|_{S_P})=3/2$.  (Remember that $P$ is not $v$-adically isolated on $S_P$, so the tangent lines are defined over the completion $k_v$.)
	\item If $S_P$ is irreducible and there is only one tangent line to $S_P$ at $P$, then $\alpha(P,L|_{S_P})=3/2$.	
\end{itemize} 

If you match up this list with the statement of the theorem, you hit paydirt: the values of $\alpha(P,L|_{S_P})$ must equal the values of $\alpha(P,L)$ when they're less than $2$, and we already knew $\alpha(P,L)\geq 2$ in all the other cases.

It therefore remains only to show $\alpha(P,L)\leq 2$.  In other words, we need to find a curve $C$ such that $\alpha(P,L|_C)\leq2$.  

It's tempting to think ``oh, I'll just use the conic I get when I intersect the plane through $P$ and $T$ with $X$''.  That pretty useful, for sure, but what if that ``conic'' is actually the union of two irrational lines?  Problems then, because that ``conic'' doesn't actually give you an upper bound for $\alpha(P,L)$.

So we roll up our sleeves and get to work.  First off, if $X$ is singular along all of $T$, then the hyperplane through $P$ and $T$ is the union of $T$ with a {\em line}, so $\alpha(P,L|_{S_P})=1$.  Thus, in what follows, we will assume that $X$ has only isolated singularities, all supported on $T$ by assumption.

Let $f\colon\tilde{X}\to X$ be a minimal desingularization of $X$, let $-K=f^*\O(1)$ be the anticanonical class on $\tilde{X}$, and let $\tilde{T}$ be the class on $\tilde{X}$ of the strict transform of the $k$-rational line $T$.  It's classical (see for example \cite{BW}) that the singularities of $X$ are canonical, so $f$ is a crepant resolution of $X$: $K_{\tilde{X}}=f^*K_X$.  This means that since $T.(-K_X)=1$, we get $\tilde{T}.(-K_{\tilde{X}})=1$ (curves contracted by $f$ intersect $K_{\tilde{X}}$ trivially), and so by adjunction we have $\tilde{T}^2=-1$.  

Let $A$ be the class $\tilde{T}-K$, so that $\deg_{-K}(A)=A.(-K)=4$ and $A^2=4$.  It is a straightforward calculation that $A$ is big and nef and that $H^0(A)=5$.  There is therefore a pencil $\mathcal{F}\cong\P^1$ of curves in $|A|$ such that each curve in $\mathcal{F}$ has a singularity at $P$.  Moreover, adjunction tells us that a smooth curve in $|A|$ has genus 1, so a general curve in $\mathcal{F}$ has geometric genus zero with a double point at $P$.  (Bertini's Theorem tells us that all but finitely many elements of $\mathcal{F}$ are irreducible.)

No two curves in $\mathcal{F}$ can have the same tangent line at $P$ because $A^2=4$, and so the morphism $\Psi\colon\mathcal{F}\to\mbox{Sym}^2(\mathcal{H})$ is injective. 

(Um.  Forgot to tell you about $\Psi$, sorry.  $\mathcal{H}$ is the projective line parametrizing tangent directions at $P$ in the Zariski tangent space at $P$.  $\Psi$ is given by $\Psi(C)=\{t_1,t_2\}$, where $t_1$ and $t_2$ represent the two tangent directions to $C$ at $P$.)

But $\mathcal{F}$ and $\mathcal{H}$ are projective curves!  So ``$\Psi$ is injective'' means ``$\Psi$ is surjective'' too.  Thus, every tangent direction at $P$ must occur as a tangent direction of some curve in $\mathcal{F}$.

We will find a curve in $\mathcal{F}$ with irrational tangent lines that are defined over $k_v$.  Theorem~\ref{thm:curve} will tell us that this curve $C$ will be the one we're looking for, with $\alpha(P,L|_C)=\deg C/2=4/2=2$.  

If $B\subset\mbox{Sym}^2(\mathcal{H})$ is the image of $\Psi$, then its preimage $\tilde{B}\subset\mathcal{H}\times\mathcal{H}$ is in the class $(1,1)$, and therefore induces an involution on $\mathcal{H}$ whose quotient is a surjective double cover $q\colon\mathcal{H}\to\mathcal{F}$.  That's a fancy way of saying ``$q$ maps (a pair of tangent directions) to (the curve in $\mathcal{F}$ that they come from)''.

All we need now is a $k$-rational element $C\in\mathcal{F}$ such that $q^{-1}(C)$ is the union of two distinct irrational points that are defined over the completion $k_v$.  In fact, we will show that there are infinitely many of these.

To do this, it is enough to prove the following lemma, which to be honest is probably not new:

\begin{lem}\label{lem:hilbert}
	Let $f\colon\P^1\to\P^1$ be a morphism of degree two, defined over a number field $k$.  For any place $v$ of $k$, there are infinitely many points $P\in\P^1$ such that $f^{-1}(P)$ consists of two points that are defined over $k_v$ but not $k$.
\end{lem}

%Since $\mathcal{F}$ and $\mathcal{H}$ both have a $v$-adically dense set of $k$-rational points, it immediately follows (by Hilbert's Irreducibility Theorem) that there is a $k$-rational element $C\in\mathcal{F}$ such that $q^{-1}(C)$ is the union of two distinct irrational points that are defined over the completion $k_v$.

\noindent
{\it Proof of lemma:} Let $Q\in\P^1(k)$ be a rational point where $f$ is not ramified.  Then there is an open $v$-adic neighbourhood $U$ of $Q$ such that for all $R\in U$, $f$ is unramified at $R$.

The set $V=f(U)$ is therefore an open $v$-adic neighbourhood of the $k$-rational point $f(Q)$ of $\P^1$, with the property that for every point $P\in V$, the set $f^{-1}(P)$ consists of two points defined over $k_v$.  The set $V$ contains a positive proportion of the $k$-rational points of $\P^1$ by height, but the Hilbert Irreducibility Theorem (see for example section 9 of \cite{Se}) implies that the proportion of $k$-rational points of $\P^1$ by height whose preimage is two $k$-rational points is zero.  Therefore there are infinitely many points of $V$ whose preimage is defined over $k_v$ but not $v$.  \qed

\vspace{.1in}

We therefore have a curve $C$ that is defined over $k$, irreducible, genus zero, singular at $P$, and has tangent directions at $P$ that are not defined over $k$, but are defined over the completion $k_v$.  (Indeed, we have infinitely many such curves.)  By Theorem~2.6 of \cite{MR2}, $C$ contains a sequence with approximation constant $\deg(C)/2=2$, and is therefore the curve we seek.  

Thus, there is always a curve of best approximation to $P$, as desired, and the value of $\alpha$ is as advertised.  \qed

\vspace{.1in}

Having treated the curve and surface cases, we will henceforth assume that the dimension of $X$ is at least three.

\begin{thm}\label{cubic}
Let $n\geq 4$ be an integer, $X\subset\P^n$ a smooth cubic hypersurface defined over $k$, $P\in X(k)$ any $k$-rational point, and let $L$ be the hyperplane class.  Assume that there is a line on $X$ defined over $k$, and let $S_P$ be the intersection of $X$ with the tangent hyperplane at $P$.  Then:
\[\alpha_L(P)=\left\{\begin{array}{ll}
1 & \mbox{if $P$ lies on a $k$-rational line} \\
2 & \mbox{if $P$ is an isolated point in $S_P(k)$} \\
3/2 & \mbox{otherwise}
\end{array}\right.\]
In all cases, there is a curve of best approximation to $P$.
\end{thm}

\begin{rem}
Note that the three cases are mutually exclusive, as any $k$-rational line through $P$ must lie on $S_P$.
\end{rem}

\noindent
{\it Proof:} First, notice that by Theorem~\ref{thm:curve}, if $P$ lies on a $k$-rational line, then $\alpha_L(P)=1$, and the line is a curve of best approximation.  Thus, we assume for the rest of the proof that $P$ does not lie on a $k$-rational line.

Let $\pi\colon Y\to X$ be the blowup of $X$ at $P$, with exceptional divisor $E$.  Then $\pi^*L-2E$ is effective (the tangent hyperplane has multiplicity 2 at $P$), with stable base locus contained in $S_P$.  By Theorem~\ref{thm:Liouville-bound}, this means that $\alpha(P,L)\geq2$, unless there is a sequence contained in $S_P$ whose approximation constant is less than $2$.  

We first consider the case where $P$ is an isolated rational point on $S_P$.  We just showed that $\alpha(P,L)\geq 2$.  All that's left is to show that $\alpha(P,L)\leq 2$.

Let $H$ be any three-dimensional linear space that intersects the tangent space at $P$ properly and transversely, contains $P$ and the line $\ell$, and is defined over $k$.  Let $V=H\cap X$ be the intersection of $H$ and $X$.  By Bertini's Theorem we may choose $H$ so that $V$ is a cubic surface in $H$ that is defined over $k$ and smooth at $P$, and whose singularities are supported on $\ell$.  Thus, by Theorem~\ref{cubicsurface}, we deduce that $\alpha(P,L)\leq\alpha(P,L|_V)\leq 2$.

We may therefore assume that $P$ is not an isolated rational point on $S_P$, and restrict our attention to $S_P$.  We will show that $\alpha(P,L|_{S_P})=3/2$, and find a curve of best approximation to $P$ with respect to $L|_{S_P}$.  (Remember we assumed that $P$ doesn't lie on a $k$-rational line!)

Let $\phi\colon S_P\dashrightarrow \P^{n-2}$ be the linear projection from the point $P$, considered inside the tangent hyperplane $T_P\cong\P^{n-1}$.  The fibres of $\phi$ (over $\kbar$) are either single points or lines. 

If $\phi$ is not dominant, therefore, then all the fibres are lines and $S_P$ is a cone.  Since $P$ is not an isolated rational point on $S_P$, there is another rational point on $S_P$, and therefore a $k$-rational line on $P$, which contradicts our earlier assumption.

Thus, we may assume that $\phi$ is dominant.  In this case, we may change coordinates in $T_P$ (which does not change $\alpha$) so that $P=[0:0:\ldots:1]$  and the defining equation for $S_P$ is
\[f(x_1,\ldots,x_{n-1})+x_ng(x_1,\ldots,x_{n-1})=0\]
where $f$ is a homogeneous cubic and $g$ is a nonzero homogeneous quadric.  The rational map
\[\psi(\x)=[x_1g(x_1,\ldots,x_{n-1}):\ldots:x_{n-1}g(x_1,\ldots,x_{n-1}):-f(x_1,\ldots,x_{n-1})]\]
is inverse to $\phi$, defined everywhere except $f=g=0$.  Thus, $S_P$ is the blowup of $\P^{n-2}$ along the scheme $f=g=0$, with the strict transform of $g=0$ contracted.  

More precisely, let $Z\subset T_P\times\P^{n-2}$ is the closure of the graph of $\phi$.  If $\pi_i$ is the projection of $Z$ onto the $i$th factor, then
\begin{itemize}
	\item $\pi_2\colon Z\to \P^{n-2}$ is the blowup along $f=g=0$.
	\item $\pi_1\colon Z\to S_P$ is an isomorphism away from the strict transform $C$ of $g=0$ from $\pi_2$ and contracts $C$ to the point $P$.
\end{itemize}

Now let $\{y_i\}\subset S_P$ be any sequence of rational points converging to $P$.  We will show that $\alpha(\{y_i\},P,L)\geq 3/2$.  The idea is to relate the line bundle $L$ to the line bundle $\O(1)$ on $\P^{n-2}$, because we know all about $\alpha$ on projective space already.

First, note that all of the $y_i$ lie in the subset of $S_P$ on which $\phi$ is an isomorphism, since none of them is $P$, and none of the lines through $P$ are $k$-rational.  That means there is a unique sequence $\{z_i\}$ of rational points on $Z$ such that $\pi_1(z_i)=y_i\in S_P$.

Now let's pull back both $L$ and $H=\O_{\P^{n-2}}(1)$ to $Z$ and compare them there.

Let $L'=L|_{S_P}$, a divisor class on $S_P$, and let $B=\pi_1^*L'$, a divisor class on $Z$.  

Let $H=\O_{\P^{n-2}}(1)$, a divisor class on $\P^{n-2}$.

Then $B=\pi_2^*H+C$ and $C=2\pi_2^*H-D$ for some effective divisor $D$.  (In case you forgot: $C$ is the strict transform of $g=0$.)  Therefore $2B-3C=2\pi_2^*H-C=D$, which is effective.  

Therefore, up to multiplication by bounded functions compactly supported away from zero, we have for any $\gamma<3/2$:
\begin{align*}
H_{L'}(y_i)\dist(P,y_i)^{\gamma} &= H_B(z_i)\dist(P,y_i)^{\gamma} \\
&\geq H_C(z_i)^{3/2}\dist(P,x_i)^{\gamma} \\
&\geq \dist(P,y_i)^{-3/2}\dist(P,y_i)^{\gamma} \\
&\rightarrow \infty
\end{align*}
as $i\to\infty$.  We conclude that $\alpha_{L'}(\{y_i\},P)\geq 3/2$, and therefore $\alpha_{L'}(P)\geq 3/2$, as desired.

All that's left is to show that $\alpha_{L'}(P)\leq 3/2$, which we can do by exhibiting an approximating sequence that achieves $\alpha=3/2$.  We can do that by finding a curve through $P$ of $L'$-degree 3 and a suitable singularity at $P$.

Since $P$ is not an isolated rational point on $S_P$, it follows that the quadric $g=0$ must have points defined over the completion $k_v$ of $k$ at the place $v$.  Let $Q$ be an irrational point of $g=0$ (which we may choose to be away from the locus $f=0$), defined over a quadratic extension of $k$ that is contained in $k_v$, and let $T$ be the line joining $Q$ to its conjugate.  

The line $T$ is defined over $k$, so let $T'$ be the closure of the curve $\psi(T)$ in $S_P$.  The map $\psi$ is an isomorphism away from $g=0$, which is contracts to $P$, so $T'$ is just $T$ with the point $Q$ glued to its conjugate.  The $L'$ degree of $T'$ is 3, so $T'$ is a cubic curve $T'$ with a node at $P$, with tangent lines corresponding to $Q$ and its conjugate.  Since the tangent lines to $T'$ at $P$ are defined over $k_v$ but not over $k$, by Theorem~\ref{thm:curve}, there is a sequence on $T'$ with $\alpha=3/2$, and $T'$ is a curve of best approximation to $P$.  \qed


\begin{thebibliography}{MR}

\bibitem[BD]{BD} Brandes, J. and Dietmann, R.; ``Rational lines on cubic hypersurfaces'', Math. Proc. Cambridge Philos. Soc. {\bf 171} (2021), no. 1, 99--112.

\bibitem[BW]{BW} Bruce, J.W.; Wall, C.T.C, ``On the classification of cubic surfaces'', J. London Math. Soc. (2) {\bf 19} (1979), no. 2, 245--256.

\bibitem[DW]{DW} Dietmann, R. and Wooley, T.; ``Pairs of cubic forms in many variables'', Acta Arith. 110 (2003), no. 2, 125–-140.

\bibitem[McK]{McK} McKinnon, David, ``A conjecture on rational approximations to rational points", J. Algebraic Geom., 16 (2007), 257--303.

\bibitem[MR]{MR} McKinnon, D.~and Roth, M., ``Seshadri constants, Diophantine approximation, and Roth's theorem
for arbitrary varieties'', {Invent. Math.} (200), 513--583 (2015). DOI 10.1007/s00222-014-0540-1.

\bibitem[MR2]{MR2} McKinnon, D. and Roth, M., ``An analogue of Liouville's Theorem and an application to cubic surfaces'', Eur. J. Math. 2 (2016), no. 4, 929--959. DOI:10.1007/s40879-016-0113-5

\bibitem[Se]{Se} Serre, J.-P., {\em Lectures on the Mordell-Weil Theorem}, Vieweg, 1997.

\end{thebibliography}
\end{document}